\documentclass[runningheads]{llncs}

\usepackage{graphicx}
\usepackage[color=pink]{todonotes}
\usepackage{booktabs}

\usepackage{hyperref}
\usepackage{amsmath}
\usepackage{xspace}
\usepackage{listings}
\usepackage[T1]{fontenc}  
\usepackage[scaled=0.80]{beramono}  
\usepackage{listings}
\lstdefinelanguage{scala}{
  alsoletter={@=>},
  morekeywords={
    abstract, case, catch, class, def, do, else, extends, false, final, finally,
    for, if, implicit, import, match, new, null, object, override, package,
    private, protected, requires, return, sealed, super, this, throw, trait, try,
    true, type, val, var, while, with, yield, domain, postcondition,
    precondition,invariant, constraint, assert, forAll,  _, return, @generator,
    ensure, require, ensuring,=>, Real, certainly, possibly, certify, errorBound,
    assertBound
  },
  sensitive=true,
  morecomment=[l]{//},
  morecomment=[s]{/*}{*/},
  commentstyle=\color{gray},
  showstringspaces=false,
  columns=fullflexible,
  mathescape=true,
  numberstyle=\tiny,
  basicstyle=\small\ttfamily,
  numbersep=5pt,
  stepnumber=2,
  numbers=none,                   
  morestring=[b]"
}
\usepackage{makecell}
\usepackage{array}
 \usepackage{multirow}
 \usepackage{rotating}
\newcommand{\libm}[0]{\texttt{libm}\xspace}

\begin{document}
\title{Sound Approximation of Programs with Elementary Functions}

\author{Eva Darulova\inst{1} \and
Anastasia Volkova\inst{2}}
\authorrunning{Darulova et al.}
%
\institute{MPI-SWS, \email{eva@mpi-sws.org}\\
 \and
Univ Lyon, Inria, CNRS, ENS de Lyon, Universit\'e Claude Bernard Lyon 1,\\ LIP UMR 5668, \email{anastasia.volkova@inria.fr}}
\maketitle              
\begin{abstract}
  Elementary function calls are a common feature in numerical programs. While
  their implementions in library functions are highly optimized, their
  computation is nonetheless very expensive compared to plain arithmetic. Full
  accuracy is, however, not always needed. Unlike arithmetic, where the
  performance difference between for example single and double precision
  floating-point arithmetic is relatively small, elementary function calls
  provide a much richer tradeoff space between accuracy and efficiency.
  Navigating this space is challenging. First, generating
  approximations of elementary function calls which are guaranteed to satisfy
  accuracy error bounds is highly nontrivial. Second, the performance of
  such approximations generally depends on several parameters which are unintuitive  
  to choose manually, especially for non-experts.

  We present a fully automated approach and tool which approximates elementary
  function calls inside small programs while guaranteeing overall user provided error
  bounds. Our tool leverages existing techniques for roundoff error computation
  and approximation of individual elementary function calls, and provides
  automated selection of many parameters. Our experiments show that significant
  efficiency improvements are possible in exchange for reduced, but
  guaranteed, accuracy.

\keywords{Elementary Functions  \and Approximation \and Synthesis \and Error analysis \and Floating-point Arithmetic.}
\end{abstract}

\section{Introduction}

Numerical programs face an inherent tradeoff between accuracy and efficiency.
For example, choosing a larger finite precision provides higher accuracy, but is
generally more costly in terms of memory and running time. Not all applications,
however, need a very high accuracy to work correctly. We would thus like to
compute the results with only as much accuracy as is needed, in order to save
resources.

Navigating this tradeoff between accuracy and efficiency is challenging. First,
estimating the accuracy, i.e. bounding errors, is non-trivial due to the complex
and discrete nature of finite-precision arithmetic which inevitably occurs in
numerical programs. Second, the space of possible implementations with different
performance characteristics is usually prohibitively large and thus cannot be
explored manually.

The need for automated tool support has been recognized previously. Today, users
can choose between different tools for analyzing accuracy of floating-point
programs~\cite{Daisy,Fluctuat,FPTaylor,Gappa,Precisa,real2Float} as well as for choosing between
different precisions~\cite{FPTuner,DaisyTuning}. The latter tools perform
mixed-precision tuning, i.e. they assign different floating-point precisions to
different operations, and thus allow to improve the performance of a program
w.r.t. a uniform precision implementation. The success of such an optimization is
usually limited to the scenario when one uniform precision is just barely not
enough to satisfy a given accuracy specification.

Another possible target for performance optimizations are elementary functions,
such as sine and exponential.
Users by default choose library function implementations which are correctly
rounded for single or double precision. Such implementations are, however,
expensive. When such high accuracy is not needed, we can save significant
resources by replacing elementary function calls by coarser approximations.
Unfortunately, existing automated approaches~\cite{StokeFloat,Metasketching} do
not provide sound accuracy guarantees.

On the other hand, tools like Metalibm~\cite{Brunie2015} provide a procedure for
approximating individual elementary functions by polynomials with rigorous
accuracy guarantees. They, however, do not consider entire programs and leave
the selection of its parameters to the user, limiting its usability mostly
to experts.
  
We present an approach and a tool which leverages an existing whole-program error
analysis and Metalibm's elementary function approximator to provide both sound
whole-program guarantees as well as efficient implementations for programs 
with elementary function calls. Given a user-provided target error
specification, our tool fully automatically distributes the error budget among
the floating-point implementation of arithmetic operations and the elementary
functions, and selects a suitable polynomial degree for Metalibm to use.

We have implemented our approach and evaluate it on several benchmarks from
literature and compare the performance of generated programs against programs
using library implementations. Naturally, we cannot compete against the highly
(often hand-optimized) library implementations for target error specifications close to 
correct-rounding. When such a high precision is not required, however, our tool 
allows users to trade performance for larger, but guaranteed, error bounds. Our
tool improves performance by on average 14\% and up to 27\% when approximating individual
elementary function calls, and on average 17\% and up to 34\% when approximating
several function calls at once. These performance improvements incur overall
whole-program errors which are only 2-3 magnitudes larger than double-precision
implementations using \libm functions and are well below the errors of single-precision
implementations.

\paragraph{Contributions}
In summary, in this paper we present
\begin{itemize}

  \item the first approximation technique for elementary functions with
    sound whole-program error guarantees,

  \item an efficient heuristic to select a suitable degree of the polynomial approximation,

  \item extensive experimental evaluation on benchmarks from literature, and

  \item an implementation, which will be released as open-source.
\end{itemize}

\section{Overview}

We will illustrate our approach on a small example from the Axbench benchmark
set~\cite{Axbench}, which computes a forward kinematics equation:
\begin{lstlisting}
def forwardk2jY(theta1: Real, theta2: Real): Real = {
  require(-3.14 <= theta1 && theta1 <= 3.14 && -3.14 <= theta2 && theta2 <= 3.14)
  val l1: Real = 0.5
  val l2: Real = 2.5

  l1 * sin(theta1) + l2 * sin(theta1 + theta2)
}
\end{lstlisting}

Although this equation is relatively simple, it still presents an opportunity
for performance savings, especially when it is called often, e.g. during the
motion of a robotics arm. Based on a rough estimate
(see~\autoref{fig:ArithOverhead}), the two elementary functions take approximately
34\% of the overall running time.

Assuming library implementations for sine, our static analysis determines the
worst-case absolute roundoff error of the result to be 3.44e-15. If such a high
accuracy is not desired, the user can specify in the postcondition a higher
error:
\begin{lstlisting}
} ensuring(res => res +/- 1e-13)
\end{lstlisting}

Our tool computes how much of the total error budget can be used to approximate
each elementary function call and calls Metalibm to generate a (different)
approximation for each elementary function call in less than 5 minutes. In this
case it determines that polynomial approximations with degrees 12 and 16, respectively, are suitable.
This new implementation is approximately 2.9\% faster than the implementation
with \libm functions. Our tool guarantees that the specified error bound is
satisfied, but in fact it often provides tighter bounds, for instance here the
actual error bound of the implementation with approximations is 2.09e-14, i.e.
roughly an order of magnitude higher than the \libm error.

This performance improvement is not very significant yet, but if we increase the
error bound further, our tool generates programs which are 13.4\% and
17.6\% faster than a \libm-based implementation (choosing degrees 24 and 20, and 
degrees 24 and 16, respectively).
The actual errors of the generated implementation are 2.21e-13 and 1.56e-12,
respectively. That is, we can save roughly half of the library function
overhead, in exchange for somewhat larger errors. These errors are still much
smaller than if we had used a single precision implementation, which incurs a
total error of 1.85e-06.

\section{Background}

Before explaining our approach in detail, we present an overview of the
necessary background. We provide substantial amount of detail on the algorithm
employed by Metalibm, as it is important for our implementation and experiments.

\subsection{Floating-point Arithmetic}

Efficiently representing real numbers on electronic computers is not a
straightforward task. In general, each real number can only be represented in
the machine with some finite precision. One of the most popular representations
of finite-precision is floating-point (FP) arithmetic. Extensive hardware and
software library support make it convenient and efficient to use.

The IEEE 754~\cite{ieee75408} standard defines several precisions and rounding
modes, of which the most commonly used ones are single and double precision with
arithmetic operations in rounding-to-nearest mode.

Any operation in finite precision is susceptible to rounding errors. When
propagated through programs, those errors can amplify and lead to (catastrophic)
output accuracy degradations.
A widely used error model~\cite{High02} of FP operations over two FP
numbers $x$ and $y$ is (in the absence of overflows and with round-to-nearest rounding mode):
\begin{align}
  \label{eqn:floats}
    x \circ_{fl} y = (x \circ y)(1 + e) + d \;\text{, }\;
    \lvert e \rvert \le \epsilon_m, \lvert d \rvert \le \delta_m
  \end{align}
  where $\circ \in {+, -, *, /}$ and $\circ_{fl}$ denotes the respective
  floating-point versions. Square root follows similarly, and unary minus
  does not introduce roundoff errors.
  The machine epsilon $\epsilon_m$ bounds the maximum relative error
  for so-called normal values. Roundoff errors of subnormal values, which provide gradual underflow,
  are expressed as an absolute error, bounded by $\delta_m$.
  $\epsilon_m = 2^{-24}, \delta_m = 2^{-150}$ and $\epsilon_m = 2^{-53}, \delta_m = 2^{-1075}$
  for single and double precision, respectively.

This model provides a convenient way to analyze small numerical expressions and
provides tight error bounds. However, manual analysis of whole programs can
quickly become tedious, and automated tool support is essential.

\subsection{Roundoff Error Analysis in Daisy}
  
  Daisy~\cite{Daisy} is a static analysis-based tool for soundly bounding
  worst-case absolute finite-precision roundoff errors:
  \begin{align*} 
    \max_{x \in [a, b]}\;\;\lvert f(x) - \tilde{f}(\tilde{x}) \rvert
  \end{align*}
  where $f$ and $x$ are a mathematical real-valued arithmetic expression and
  variable, respectively, and $\tilde{f}$ and $\tilde{x}$ are their finite-precision
  counterparts. This definition extends to multivariate $f$ component-wise.

  Worst-case roundoff errors substantially depend on the ranges of inputs which
  have to be provided by a user ($[a, b]$ is the range for $x$ given in such a
  precondition in the equation above). Roundoff errors due to arithmetic
  operations equally depend on the ranges, but these cannot be easily provided
  by the user.

  Given an arithmetic expression, together with an input domain, Daisy performs
  a two-step dataflow static analysis. In the first step, it uses
  interval~\cite{intervals} or affine arithmetic~\cite{affineArithmetic} to
  compute the ranges of all intermediate expressions, and in the second step it
  uses this information to compute and propagate roundoff errors. The
  propagation is performed with affine arithmetic~\cite{affineArithmetic}.

  This analysis is efficient and computes tight error bounds for straight-line
  arithmetic expressions over $+, -, *, /, \sqrt{}$ as well as the most commonly
  used elementary function calls (sin, cos, exp, log, tan, assuming a slightly
  larger error bound).
  Daisy currently does not support loops or conditional branches; for a discussion
  on the challenges involved, we refer the reader to~\cite{Rosa,Goubault2013}.

\subsection{Standard \libm and its Limitations}

An inherent part of scientific and financial computations are mathematical functions. 
The standard mathematical library containing elementary (exp, log, sin, cos and their inverses, etc) and special (power $x^y$, erf, $\Gamma$, etc) functions is \libm. 
This library is fully specified in the C language standard (ISO/IEC 9899:2011) though the recent 2008 revision of the IEEE 754 floating-point standard has also attempted standardization.
There are various different implementations of \libm  that depend on the operating system and programming language.
Here when referring to \libm we mean the GNU \texttt{libc}\footnote{\url{https://www.gnu.org/software/libc/}} implementation.
GNU \libm provides a reference implementation of basic mathematical functions targeting double and single precision inputs. 

In scientific computing the performance of \libm is of crucial importance. 
For example, large-scale simulation codes in high-energy physics run at CERN~\cite{Innocente2012} rely on algorithms for reconstruction of particle collision coordinates that spend more than half of their time on computing the function $\exp$.
Similarly, the SPICE simulator~\cite{SPICE} of analog circuits is based on solutions of non-linear differential equations and spends up to three quarters of its time on computing elementary functions.
Finally, Figure~\ref{fig:ArithOverhead} illustrates a rough estimate of the time spend on the calls to \libm in our set of benchmarks. 

\begin{table}[t]
  \centering
\begin{tabular}{lcccccc}
\toprule
                        & \multicolumn{3}{c}{64 bit}     & \multicolumn{3}{c}{32 bit}     \\
benchmark               & \libm  & arith. only & overhead & \libm  & arith. only & overhead \\
\midrule
xu1                     & 553.5 & 367.8       & 0.34     & 449.1 & 364.6       & 0.19     \\
xu2                     & 572.6 & 313.9       & 0.45     & 407.8 & 309.7       & 0.24     \\
integrate18257          & 527.4 & 321.5       & 0.39     & 428.4 & 328.2       & 0.23     \\
integStoutemyer & 503.7 & 286.5       & 0.43     & 437.0 & 359.7       & 0.18     \\
sinxx10                 & 395.0 & 354.7       & 0.10     & 380.8 & 320.4       & 0.16     \\
axisRotationX           & 479.9 & 413.4       & 0.14     & 375.8 & 367.9       & 0.02     \\
axisRotationY           & 481.7 & 324.6       & 0.33     & 376.3 & 336.1       & 0.11     \\
rodriguesRotation       & 610.0 & 428.7       & 0.30     & 478.2 & 447.5       & 0.06     \\
pendulum2               & 583.4 & 343.8       & 0.41     & 491.3 & 317.2       & 0.35     \\
forwardk2jY             & 468.7 & 310.3       & 0.34     & 425.2 & 349.7       & 0.18     \\
forwardk2jX             & 492.0 & 310.3       & 0.37     & 398.4 & 357.2       & 0.10    \\
\bottomrule
\end{tabular}
\caption{Average number of cycles needed for executing each benchmark with and without elementary functions (for double
and single floating-point arithmetic implementations)}\label{fig:ArithOverhead}
\end{table}

Using the standard mathematical library in practice has several limitations.
Standard \libm functions were implemented to work with arguments on the whole range of representable numbers for the target format, and 
produce results which are as accurate as the target format allows. 
In practice, this is often more than is necessary.
For example, the calls to cosine in some CERN simulations~\cite{Piparo2012} always remain
within one period and, due to the noise in measurements, the result is only needed to be correct to a few digits. 
There is obviously a need for more flavors of elementary functions: using non-standard domains and target accuracies.

\subsection{Elementary Function Approximation in Metalibm}
Adapting elementary function codes for different domains, accuracy, performance metrics (e.g. throughput vs. latency) is highly non-trivial and error prone.
%
The state of the art solution for the automation of \libm development is the Metalibm tool\footnote{\url{http://www.metalibm.org}}~\cite{Brunie2015}. 
It is able to quickly generate code for elementary functions and guarantee that the final accuracy is bounded by a user-given value.
The produced code obviously cannot compete with handwritten manually-optimized codes of standard \libm (for target errors close to correct rounding).
However, the approximations are available quickly and automatically.


Metalibm is a pushbutton tool that tries to generate code, on a given domain and for a given accuracy, for evaluation of an arbitrary univariate function with continuous derivatives up to some order. 
There is no dictionary with a fixed set of functions, the target function is specified as a mathematical expression instead. 

Metalibm currently handles exponential, logarithm, cosine, sine and tangent and their inverses, and hyperbolic cosine, sine and tangents, as well as their inverses.  
There has been some work towards support of special functions, such as Bessel~\cite{Lauter2015}. 
Unfortunately, the bivariate functions are out of reach at the moment.

As backend, Metalibm uses the Sollya scripting language/tool~\cite{ChevillardJoldesLauter2010}.
Among others, Sollya provides state-of-the-art polynomial approximations~\cite{FPminimax}.
To guarantee the bounds on evaluation errors, Metalibm uses the Gappa proof assistant\footnote{\url{http://gappa.gforge.inria.fr}}~\cite{Gappa2011}
which uses a very similar error analysis as Daisy.

\subsubsection*{Parameter space}

The goal of Metalibm is to provide more choices in implementations of metamathematical functions.
The most important choices include:
\begin{itemize}
	\item function $f$ itself
	\item domain of implementation $[a;b]$
	\item requirements on the accuracy of the final code
	\item maximum degree of the approximation polynomial $d_\text{max}$  
	\item size $t$ (bits) of the table index for table-driven methods
\end{itemize}

However, Metalibm does not help to navigate through this parameter space, proposing no guidelines to the user.
This limits the usage of the tool mostly to experts in \libm development.

\subsubsection*{High-level Algorithm}
Metalibm generates approximations in the following steps:
\begin{itemize}
	\item \emph{argument reduction:} First, Metalibm tries to detect algebraic properties of functions, such as parity or periodicity, in order to exploit them for an argument reduction. 

	\item \emph{domain splitting:} It may happen that even on possibly reduced argument approximation with one polynomial of maximum degree $d_\text{max}$ is impossible. 
	Then, Metalibm splits the domains and performs a piece-wise polynomial approximation. 
	 
	 \item \emph{polynomial approximation:} Finally, Metalibm computes the polynomial approximation(s) for the reduced domain (or subdomains for splitting) and generates Gappa proofs to account for the approximation and evaluation errors due to floating-point arithmetic. Polynomials are evaluated using Horner's scheme.
	
\end{itemize}
At the end, Metalibm produces efficient C code implementing the target function. 

\subsubsection*{Argument reduction and properties detection}

Exploiting algebraic properties helps to reduce the approximation domain in order to decrease the minimum degree of the approximation polynomial. 
For simple functions it is often possible to reduce the domain such that only one polynomial is used for the approximation.
Only if the implementation with one polynomial is impossible, the tool moves to domain splitting.
In Metalibm, domain reduction is always prioritized over domain splitting: in case of splitting, we must store the splitting points and coefficients of approximation polynomials for each domain. 
In addition to that, some splitting techniques yield \texttt{if-else} statements, which decreases throughput for vectorized implementations.

The tricky part is that Metalibm does not know which function it implements analytically, so using text-book argument reduction schemes, e.g. for exponential, is impossible.
The idea is to verify properties numerically, up to some accuracy.
Thus, the domain should be ``large enough'' to be enable detection.
Currently, the following properties are supported: 
\begin{align*}
f(x+y) &= f(x)f(y) \quad &\text{exponential functions}\\
f(x+C) &= f(x)		\quad &\text{periodic functions} \\
f(x) + f(y) &= f(xy) \quad &\text{logarithmic functions} \\
f(-x) = f(x); &f(-x) = -f(x) \quad &\text{odd/even functions}			
\end{align*}
Obviously, some functions, for example some compound functions,  simply do not have efficient argument reduction schemes.
In this case, piece-wise polynomial approximation is performed.

\subsubsection*{Domain splitting}
Given a function $f$, domain $[a;b]$ and maximum polynomial approximation degree $d_{\text{max}}$, Metalibm tries to split the domain into non-overlapping intervals such that on each of them an approximation of low degree $d \leq d_\text{max}$ is possible. 
There exist different schemes for domain splitting: uniform, logarithmic or arbitrary.
Metalibm uses a non-uniform splitting~\cite{Kupriianova2014} . 
The idea is to perform splitting only when it is impossible to approximate with the maximum degree.
The algorithm uses a heuristic search for the best splitting and is based on the de la Vall\'ee Poussin theorem to compute the minimal approximation degree.
Their method results in (almost) uniform degree-usage for all subdomains, which yields uniform memory usage for polynomial coefficients and stable performance. 
However, this strategy means that Metalibm uses the polynomial approximation of the maximum degree basically all the time, and rarely less.
This leaves it to the user to select a suitable degree, mostly by trial and error. 
In addition to that, even slight changes to the implementation domain might lead to completely different domain splitting scheme and different degrees of approximations.

\subsubsection*{Polynomial approximation}

When the domains are reduced, Metalibm uses Sollya-generated Remez-like approximation polynomials.
These polynomials are guaranteed to have the best floating-point coefficients for a given precision~\cite{FPminimax}. 
Polynomial evaluation is implemented using Horner's scheme and then a Gappa proof~\cite{Gappa2011} is automatically generated to obtain the bound on the approximation error.
It should be noted that Metalibm is quite conservative in its implementations and usually leaves an error threshold between the actually implementation and user-given error bound. 

Sometimes, to improve the efficiency of approximation, table-driven methods~\cite{Muller2016} are used. The idea is to store the values of the function at some points of the domain in a table and use the function's algebraic properties to reduce the degree of the polynomial approximation. 
For each function, the choices on the usage of the tables and their size highly depend on the properties of the underlying hardware (size of caches, speed of arithmetic vs. memory access, etc.).

\subsubsection*{Reconstruction}

The goal of reconstruction is to provide mapping between the value of function $f$ on the reduced argument or a subdomain to its value on the initial larger domain.
When an argument was reduced using algebraic properties, the reconstruction is just an application of the inverse transformation.
In case of domain splitting, reconstruction is more complicated. 
The common way for arbitrary splitting is just a series of $\text{if-else}$ statements. 
However, this approach prevents automatic vectorization of generated code and Metalibm proposes a more versatile solution.
It computes a special polynomial which, when evaluated on the function's argument, gives the index of the subdomain to be used.
While being elegant, for scalar implementations this approach yields additional polynomial evaluation.
This makes the dependency between the performance the code, function domain and approximation degree highly nonlinear. 

Consider, for example, a case when on a domain $[a;b]$ with a degree $d_\text{max}$ there was no need in domain splitting. 
It may happen that slight narrowing of this domain requires a domain splitting for the same $d_\text{max}$. The implementation cost increases by at least a cost of polynomial evaluation for the reconstruction. On the other hand, increasing $d_\text{max}$ by one could have maintained implementation without domain splitting by adding a cost of just one multiplication and addition.

\section{Whole Program Approximation}

In this section we describe our approach in detail. Our tool supports
approximations of straight-line expressions featuring the standard arithmetic
operators ($=, -, *, /$) as well as a set of the most commonly used elementary
functions ($\sin{}, \cos{}, \tan{}, \log{}, \exp{}, \sqrt{}$), and attempts to approximate
the latter.
To avoid confusion, we will use the term `program' for the entire expression in
which the user would like to approximate elementary function calls, and
`function' for individual elementary functions.
In addition to specifying the expression itself, the user of our tool also
specifies the domains of all its inputs, together with a target overall absolute
error which should be satisfied.

\subsection{Overall Structure}
\begin{figure}[t]
  \includegraphics[width=\textwidth]{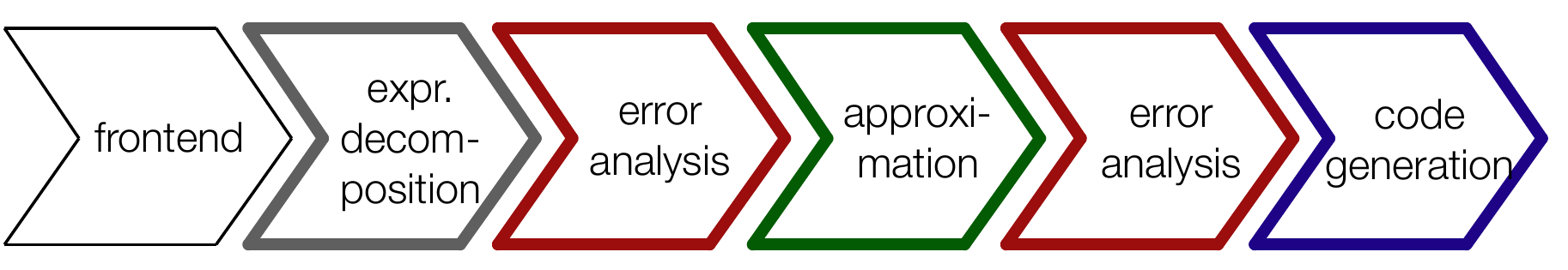}
  \caption{Overview of different steps in our approach.}
  \label{fig:pipeline}
\end{figure}

We implement our approach inside the Daisy framework~\cite{Daisy}. Daisy is built up in
phases which provide different functionalities. We re-use the phases related to
frontend and backend handling and roundoff error analysis and add two new ones
to support automated elementary function approximation. \autoref{fig:pipeline}
shows an overview of the overall structure.

After parsing the input file, Daisy decomposes the abstract syntax tree (AST) of
the program we want to approximate such that each elementary function call is
assigned to a fresh local variable. This transformation eases the later
replacement of the elementary functions with an approximation.

As the next step, Daisy runs a roundoff error analysis on the entire program,
assuming a \libm implementation of elementary functions. This analysis computes
a real-valued range, and a worst-case absolute roundoff error bound for each
subexpression in the AST. In particular, it computes the ranges and errors of
the arguments of elementary function calls, which we will use for approximation
in the step.

Now we can perform approximation of elementary functions. For this, Daisy calls
Metalibm for each elementary function which was assigned to a local variable. If
Metalibm successfully computes an approximation, it generates a C implementation.
Daisy extracts relevant information about the generated approximation (name of
file, name of function, etc.) and stores it in the AST. We discuss the details
of calling Metalibm in the following section.

Next, Daisy performs another roundoff error analysis, this time taking into
account the new approximation's precise roundoff error bound reported by
Metalibm.
Finally, Daisy generates C code for the program itself, as well as all necessary
headers to link with the approximation generated by Metalibm.

\subsection{Approximation Phase}\label{sec:approx-phase}

Several Metalibm parameters determine the accuracy and efficiency of the
generated elementary function approximations. We next discuss how Daisy handles
the selection of the degree of the polynomial, and the target error of each
elementary function call, which are two of the most important parameters.

\paragraph{Error Distribution}
  The user of our tool specifies only the error bound of the final result of each
  program, i.e. he or she does not have to specify the target error of each
  individual elementary function call separately. This is important for usability,
  as it is usually nontrivial even for experts to determine which error an
  individual elementary function should satisfy. This is especially true if there are
  several elementary function calls inside one program, or a function call is not
  the last statement of a program. Furthermore, Daisy (and all other similar
  tools) measures roundoff in terms of absolute errors, whereas Metalibm takes as
  input specification relative errors, hence a translation is needed.

  Thus, Daisy needs to distribute a total absolute error specified by the user
  among the potentially several elementary function calls, and the remaining
  arithmetic operations. Let us denote the total absolute error by
  $\tau = |f(x) - \tilde{f}(\tilde{x})|$, where $f$ denotes the real-valued
  specification of the program, and $\tilde{f}$ the final finite-precision
  implementation with elementary function approximations. We will denote the
  real-valued implementation of the program but with elementary function calls
  replaced by polynomial approximations by $\hat{f}$.

  Daisy conceptually decomposes the total error budget into an approximation and a roundoff
  error part:
  \begin{equation}\label{eqn:error-distrib}
  |f(x) - \tilde{f}(\tilde{x})| \le |f(x) - \hat{f}(x)| + |\hat{f}(x) - \tilde{f}(\tilde{x})|
  \end{equation}
  The first part ($\tau_{approx} = |f(x) - \hat{f}(x)|$) captures the error due to the elementary function
  approximations, assuming the arithmetic part of the program is still real-valued, i.e.
  without roundoff errors. The second part ($\tau_{fl} = |\hat{f}(x) - \tilde{f}(\tilde{x})|$) captures
  the roundoff error due to the finite-precision arithmetic. This separation is needed
  for Daisy to determine how much of the error can be used for the approximations.

  Note, however, that at this point, Daisy cannot compute $|\hat{f}(x) -
  \tilde{f}(\tilde{x})|$ exactly, as the approximations are not available yet.
  Instead, it computes an estimate of the roundoff error, by performing the
  analysis assuming \libm library function implementations. This error will, in
  general, be slightly smaller than the final roundoff error. On the other hand,
  Metalibm usually generates approximations which satisfy error bounds which are
  smaller than what was asked for. For this reason, Daisy performs a second round
  of error analysis, after the approximations have been determined. In this second
  round, all errors (smaller or larger) are correctly taken into account.

  If the program under consideration contains several elementary function calls,
  the approximation error budget needs to be distributed among them.
  For this, we consider two distribution strategies:
  \begin{itemize}
    \item
      With the \emph{equal} distribution strategy, Daisy distributes the error
      equally among elementary function calls. That is, for $n$ function calls,
      each will have an error budget of $\tau_i = \tau_{approx}/n$ for each call.

    \item With the \emph{derivative} based strategy, we take into account that
      that the errors introduced by elementary approximations may be propagated
      differently through the remaining program: some errors may be magnified
      more than others. It thus seems reasonable to allocate more of the error
      budget to those function calls, whose errors are potentially magnified
      more. We estimate the error propagation by symbolically computing the
      derivative w.r.t. to the elementary function call, and by bounding it
      w.r.t. to the variable ranges specified.

      The derivative bounds for the different function calls are normalized to one,
      such that we compute a weight $w_i$ for each call. The error budget for each call
      is then computed as $\tau_i = \tau_{approx} * w_i$.
  \end{itemize}

\paragraph{Total Error vs. Local Error}

  Consider a program with two elementary function calls. Expanding
  \autoref{eqn:error-distrib} we obtain
  \begin{equation*}
    |f(x) - \tilde{f}(x)| = |f(x) - \hat{f}_1(x)| + |\hat{f}_1(x) - \hat{f}_2(x)| + |\hat{f}_2(x) - \tilde{f}(\tilde{x})|
  \end{equation*}
  where $\hat{f}_1$ denotes the program with one elementary call replaced by an
  approximation and $\hat{f}_2$ with both calls replaced. The error distribution
  assigns error budgets $\tau_1$ and $\tau_2$ to the first and second term, i.e.
  to the first and second elementary function call, respectively.

  These error budgets represent the total error due to the elementary function
  call \emph{at the end of the program}. For calling Metalibm, however, we need
  the \emph{local} error at the function call site. Again, due to error
  propagation, these two errors can differ significantly, and may lead to
  overall errors which exceed the error bound specified by the user.

  We solve this issue by computing an approximation of the local error from the
  global error budget $\tau_i$. The local error gets propagated through the
  remaining program, and can thus be magnified or diminished. We compute a
  linear approximation of this propagation error $m_i$, and then compute the
  local error $\varepsilon_i$ as
  \begin{equation}
    \varepsilon_i = \tau_i / m_i   \quad \quad\text{ where } m_i = \max_{[a, b]}\frac{\partial{f}}{\partial{x_{i}}}
  \end{equation}
  where $x_i$ is the local variable to which the elementary function call is
  assigned to in the decomposition phase. Daisy computes partial derivatives 
  symbolically and maximizes them over the specified input domain.

  We proceed similarly for all elementary functions in the program. We process
  the functions in the order in which they appear in the program. This is
  important, as we want to compute the partial derivatives on the real-valued
  programs (it is not clear how to compute derivatives of approximations). We
  compute the derivatives on the program which follows the function call. Thus,
  when processing the first call, all remaining (and following) ones are still
  real-valued and the derivative computation is possible. For the next call, the
  previously introduced approximation is preceding this call, thus the
  derivative computation can still be performed over the reals.

\paragraph{Degree Selection}
  The target error is one important parameter which significantly influences the
  efficiency of the approximations generated by Metalibm. Another one is the
  degree of the approximation polynomial. Unfortunately, it is not easy to
  predict exactly which degree is optimal, as the efficiency nonlinearly depends also on
  the domain of the function call, the function which is to be approximated
  itself, the internal logic of Metalibm, etc.

  Daisy thus performs a linear search, starting from a small degree. The search stops, either when
  the (coarsely) estimated running time reported by Metalibm is significantly higher than
  the current best, or when Metalibm times out. The approximation with the smallest
  estimated running time is returned.

  We have empirically observed that for functions which are monotone (i.e.
  $\exp{}$, $\log{}$, $\sqrt{}$, $\tan{}$), small degrees often result in efficient
  approximations, whereas for non-monotone functions (i.e. $\sin$, $\cos$),
  higher degrees are necessary. For the former class, we thus search the degrees
  4, 8, 12, and 16, and if the function to be approximated contains $\sin$ or
  $\cos$, then we choose degrees 12, 16, 20 and 24. We have observed that steps
  of four are a good compromise between efficiency of implementations and
  efficiency of the analysis.

\paragraph{Depth of Approximation}

  Another parameter with potentially high impact is what we call the depth of
  the approximation. We can ask Metalibm to approximate individual elementary
  functions, but it can also approximate \emph{compound} functions (e.g. $\sin(\cos(x)-1)$) directly.
  If successful, this has the advantage that only one
  polynomial approximation is generated and only one function call is involved.
  In addition to that, the Metalibm-generated code often wins in accuracy, over standard consecutive calls to the \libm, thus provides a larger buffer to be traded for performance.  Finally, Metalibm can overcome undesired effects, e.g. cancellation at zero for the $\cos(x)-1$ above. 

  This is a highly application specific parameter, and we thus leave it
  under user control. Daisy provides automation in the sense that the user only
  needs to specify the desired depth as a number, where depth means height of the
  abstract syntax tree.

\paragraph{Other Parameters}

Metalibm provides several other parameters, such as width of the table index for table-driven methods, the minimum width of the subdomains for splitting, etc. 
Our current set of benchmarks uses the standard set of elementary functions
(exp, log, sqrt, sin, cos, tan), which are relatively ``well-behaved''. We thus
leave other parameters with their default values, which seem to perform well for
our set of functions.
For more demanding functions~\cite{Lauter2015,Brunie2015}, fine-grained manual
tuning by experts is required; automation of such specialized functions is out
of scope of this work.

\subsection{Code Generation}
Metalibm generates C code for the
elementary function approximations, and Daisy generates code for the arithmetic
part of the program. These two need to be linked for a usable implementation.
Daisy generates the necessary headers, as well as a compilation script for this
purpose. This script inlines all approximations, as we have observed that this
is necessary for performance. \libm calls also get inlined with the \texttt{-O2}
compiler flag which we use for compilation.



\section{Experimental Evaluation}

We evaluate our approach on a number of benchmarks from the literature in terms
of accuracy and performance.

\subsection{Experimental Setup}

All experiments have been performed on a Debian Linux 9 Desktop machine with a
3.3 GHz Intel i5 processor with 16 GB of RAM. We use the GNU g++ compiler
version 6.3.0. All code for benchmarking is compiled with the \texttt{-02} flag.

\paragraph{Benchmarks}
\begin{table}[t]
  \centering
  \begin{tabular}{l|>{\centering\arraybackslash}p{1.5cm} >{\centering\arraybackslash}p{1.5cm} >{\centering\arraybackslash}p{1.5cm} >{\centering\arraybackslash}p{1.5cm} >{\centering\arraybackslash}p{2cm}}
    \toprule
    \multirow{2}{*}{benchmark} & \multicolumn{1}{c}{\multirow{2}{*}{libm}} & \multicolumn{1}{c}{small} & \multicolumn{1}{c}{middle} & \multicolumn{2}{c}{large}                             \\ \cline{3-6} 
    & \multicolumn{1}{c}{}                      & \multicolumn{1}{c}{equal} & \multicolumn{1}{c}{equal} & \multicolumn{1}{c}{equal} & \multicolumn{1}{c}{deriv} \\ \hline
    sinxx10                    & 2.56e-13                                  & 2.64e-12                  & 2.51e-11                   & 2.50e-10                  & 2.50e-10                  \\
    xu1                        & 6.86e-15                                  & 1.65e-14                  & 1.62e-13                   & 2.01e-12                  & 2.38e-12                  \\
    xu2                        & 9.44e-15                                  & 1.50e-14                  & 4.21e-14                   & 3.47e-13                  & 3.79e-13                  \\
    integrate18257             & 3.68e-15                                  & 1.57e-14                  & 1.65e-13                   & 1.80e-12                  & 1.85e-12                  \\
    integStoutemyer           & 1.19e-15                                  & 3.12e-15                  & 2.67e-14                   & 2.61e-13                  & 7.91e-14                  \\
    axisRotationX              & 5.77e-15                                  & 2.51e-14                  & 1.91e-13                   & 2.19e-12                  & 1.88e-12                  \\
    axisRotationY              & 1.11e-14                                  & 1.94e-13                  & 1.88e-12                   & 2.19e-11                  & 2.50e-11                  \\
    rodriguesRotation          & 3.93e-13                                  & 2.16e-12                  & 1.90e-11                   & 2.19e-10                  & 1.95e-10                  \\
    pendulum1                  & 4.61e-16                                  & 1.91e-14                  & 1.88e-13                   & 1.88e-12                  & 1.88e-12                  \\
    pendulum2                  & 7.00e-15                                  & 1.90e-14                  & 1.53e-13                   & 1.49e-12                  & 8.48e-14                  \\
    forwardk2jX                & 7.22e-15                                  & 2.62e-14                  & 2.23e-13                   & 2.19e-12                  & 1.93e-12                  \\
    forwardk2jY                & 3.44e-15                                  & 2.09e-14                  & 2.21e-13                   & 1.56e-12                  & 1.98e-12                  \\
    ex2\_1                     & 1.28e-15                                  & 1.34e-14                  & 2.50e-13                   & 1.76e-12                  & 1.76e-12                  \\
    ex2\_2                     & 2.66e-15                                  & 2.28e-14                  & 2.20e-13                   & 2.50e-12                  & 2.50e-12                  \\
    ex2\_3                     & 1.33e-15                                  & -                         & 2.50e-13                          & 2.50e-12                  & 2.50e-12                  \\
    ex2\_4                     & 2.11e-15                                  & 1.38e-14                  & 1.14e-13                   & 1.76e-12                  & 2.50e-12                  \\
    ex2\_5                     & 1.09e-14                                  & 8.93e-14                  & 5.39e-13                   & 8.06e-12                  & 8.65e-12                  \\
    ex2\_9                     & 3.37e-15                                  & 2.05e-14                  & 1.89e-13                   & 1.88e-12                  & 1.88e-12                  \\
    ex2\_10                    & 3.00e-15                                  & 1.79e-14                  & 1.59e-13                   & 1.56e-12                  & 2.50e-12                  \\
    ex2\_11                    & 6.08e-14                                  & 9.36e-14                  & 6.04e-13                   & 5.70e-12                  & 5.12e-12                  \\
    ex3\_d                     & 2.78e-15                                  & 9.78e-15                  & 9.79e-14                   & 1.29e-12                  & 1.29e-12                  \\ \bottomrule
  \end{tabular}
  \caption{Accuracy comparison of approximated programs against libm.}
  \label{tbl:accuracy}
\end{table}

The benchmarks forwardk2j* are part of the Axbench approximate computing
benchmark suite~\cite{Axbench}. Inspired by this benchmark, we have included axisRotation*,
which rotates the x-y axis in carthesian coordinates counter-clockwise by theta
degrees. RodriguesRotation is a more involved formula for rotating a vector in
space\footnote{\url{https://en.wikipedia.org/wiki/Rodrigues27_rotation_formula}}. 
The pendulum benchmarks have been previously used in the context of
roundoff error verification~\cite{Rosa}. The benchmarks sinxx10, xu1, xu2, integrate18257,
integrateStoutemyer are from the benchmark set of the COPRIN project
\footnote{\url{http://www-sop.inria.fr/coprin/logiciels/ALIAS/Benches/}}.
The benchmarks ex2* and ex3\_d are problems from a graduate analysis textbook.
We show all benchmarks used in the appendix; while they are relatively short,
then represent important kernels usually employing several elementary function calls.

We have based the target error bounds on the roundoff errors obtained for an
implementation which uses \libm library functions. 
For each benchmark we have three sets of target errors: small, middle and large errors, each of which is roughly two, three and four orders of magnitudes larger than the \libm-based bound, respectively. 
\autoref{tbl:accuracy} shows the error bounds of the \libm-based implementations, 
as well as the final error bounds obtained by the implementation with
approximations. Note that the table does not show the target error bounds, but
the actually achieved ones. The exact target error bounds are listed in the appendix.

\paragraph{Performance Benchmarking}
To provide a fine-grained comparison, we measure performance in terms of processor clock cycles. 
We use \texttt{RDTSC} instruction available on x86 architectures. This function returns a time-stamp counter increased every clock cycle. We ensure that the instruction pipeline is flushed to avoid out-of-order executions and perform cache warm-up. 

Daisy provides all necessary infrastructure for test case code generation and
automatic benchmarking. Each benchmarking executable runs the Daisy-generated
code on 100000 random inputs from the input domain. Of the measures number of
cycles we discard the highest 10\%, as we have observed these to be outliers.
For the remaining ones, we compute the average, minimum and maximum number of
cycles Finally, we perform 5 benchmarking runs for each test case and report
average results of those.

\subsection{Experimental Results}

\begin{table}[t]
\centering
\begin{tabular}{lrr|rrrrr|r}
\toprule
                        & small    & middle  & \multicolumn{5}{c|}{large}                                & \multirow{2}{*}{float 32} \\
benchmark               & equal    & equal   & equal   & deriv   & no table & depth 1  & depth $\infty$ &                           \\
\midrule
axisRotationX           & 12.6   & 18.9  & 10.0  & 18.0  & 11.1   & -      & -            & -10.1                   \\
axisRotationY           & 5.0    & 16.7  & 14.0  & 18.9  & 13.4   & -      & -            & -16.0                   \\
ex2\_1                  & 9.9    & 12.4  & 13.2  & 13.1  & 11.8   & 12.8   & 11.5         & 7.3                     \\
ex2\_10                 & 22.2   & 21.7  & 22.2  & 9.7   & 22.2   & 21.7   & 21.9         & 7.2                     \\
ex2\_11                 & -13.5  & -6.3  & -2.5  & -4.9  & -0.9   & -4.4   & -6.1         & 18.6                    \\
ex2\_2                  & 0.2    & 9.1   & 18.2  & 18.7  & 17.7   & 19.2   & 18.9         & 4.5                     \\
ex2\_3                  & -      & 26.7  & 27.4  & 26.8  & 27.3   & 26.8   & 27.4         & 5.7                     \\
ex2\_4                  & 11.0   & 14.0  & 14.8  & 13.8  & 14.8   & 14.9   & 13.5         & 7.0                     \\
ex2\_5                  & 24.3   & 24.1  & 25.4  & 25.0  & 24.5   & 25.6   & 26.1         & 5.3                     \\
ex2\_9                  & 2.2    & 10.2  & 8.2   & 9.5   & 9.9    & 8.0    & 10.0         & -8.7                    \\
ex3\_d                  & 4.7    & 9.4   & 20.0  & 20.6  & -0.3   & 20.4   & 20.6         & 5.8                     \\
forwardk2jX             & -0.6   & 17.1  & 17.6  & 18.4  & 17.7   & -      & -            & -8.3                    \\
forwardk2jY             & 2.9    & 13.4  & 17.6  & 17.9  & 15.3   & -      & -            & 10.2                    \\
integrate18257          & -10.6  & 15.6  & 22.3  & 13.4  & -3.7   & 18.1   & 34.1         & 0.8                     \\
integStoutemyer & -11.8  & 3.8   & 6.0   & 1.3   & -6.0   & 4.4    & 19.8         & 6.7                     \\
pendulum1               & -6.0   & -4.5  & -4.5  & -4.2  & -7.8   & -      & -            & -9.2                    \\
pendulum2               & 9.1    & 9.6   & 11.5  & 5.7   & -0.0   & 21.2   & 19.7         & 1.4                     \\
rodriguesRotation       & 9.6    & 16.1  & 14.7  & 13.9  & 14.6   & 10.8   & 11.6         & -6.6                    \\
sinxx10                 & 6.4    & 8.3   & 8.1   & 8.3   & 8.3    & -15.5  & 1.3          & 6.3                     \\
xu1                     & -22.5  & 14.7  & 26.7  & 24.5  & 27.5   & 27.4   & 26.3         & 8.5                     \\
xu2                     & -1.5   & 3.5   & 11.7  & 12.3  & 12.1   & 11.3   & 25.2         & -3.2                    \\
\midrule
average                 & 2.7    & 12.1  & 14.4  & 13.4  & 10.9   & 13.9   & 17.6         & 1.6                    \\
average success         &9.2 & 13.9 &  16.2 &  15.2 &  16.2&  17.1 &  19.1 &  6.5 \\
\bottomrule
\end{tabular}
\caption{Performance improvements (in percent) of approximated code w.r.t. a program with \libm library function calls.}
  \label{tab:performance}
\end{table}

\paragraph*{Performance Improvements}
\autoref{tab:performance} shows the performance improvements of approximated code w.r.t. \libm based implementations of our benchmarks.
We show improvements for the following settings:
\begin{itemize}
  \item for small, middle and large target errors with an equal error distribution, and the table index width set to 8 bits (which enables table-based methods), approximating each elementary function call separately
  \item for large target errors with a derivative-based error distribution (and table index width 8) and individual approximation
  \item for large target errors with an equal error distribution, but with no table (table index width set to 0), with individual approximation
  \item for large target errors, approximation of compound functions with depth 1
  \item for large target errors, approximation of compound functions with as much depth as is possible (until the expression is no longer univariate)
\end{itemize}
For benchmark `ex2\_3', Metalibm timed out for all elementary function calls. For the compound approximations, we only report 
performance improvements for functions which actually have compound calls.

We observe that our tool can generate code with significant performance
improvements for most functions. As expected, the improvements are smallest for the
tightest target error bounds, and largest for the largest ones, on average 2.7\%
to 14.4\%. If we consider only those cases where the approximation was successful,
i.e. the performance is better than a \libm implementation (otherwise we can just use
\libm), the average improvements are between 9.2\% and 16.2\%.
We further note that for several benchmarks, improvements of more than 20\% (up
to 27\%) are possible. Comparing the performance improvements to our estimates
of the portion of running time spent in library calls
(see~\autoref{fig:ArithOverhead}), we see that we can recuperate much of the
overhead (in exchange for lower accuracy, see~\autoref{tbl:accuracy}).

Somewhat surprising, we did not observe an advantage of using the
derivative-based error distribution over the equal one. We suspect that is due
to the nonlinear nature of Metalibm's heuristics. We thus leave the equal error
distribution as the default one.

\autoref{tab:performance} further demonstrates the usage of tables generally improves the performance. 
Only for the benchmark xu1 we observe a slightly better performance when not using a table.
Furthermore, we observed in our tests that for our magnitudes of target errors
increasing the table sizes often leads to slower implementations (not shown in
the table): tables then do not fit in L1 or L2 cache and memory access time
prevails comparing to the computation time.

We can clearly observe the nonlinear relation between the performance and the target error for the approximations (for equal distributions they grow linearly with the small/middle/large program errors).
For instance, for axisRotationX and rodriguesRotation benchmarks increasing the target errors from low to middle provides clear improvement but passing from middle to large leads to a worse performance. 
This is the result of discrete decisions concerning the approximation degrees and the domain splittings inside Metalibm. 

For those benchmarks which contain compound functions, we observe that it is
generally beneficial to approximate `as much as possible'. However, for some
benchmarks, like sinxx10, we observe that Metalibm cannot provide any
improvement. In general, we noticed that it has trouble in approximating sine
functions and its powers.

Finally, we also considered implementations of our benchmark set in single floating-point precision,
even though Metalibm does not support implementations in single precision (but can target large errors).
We considered the following setting: program target errors are scaled to single precision, inputs/outputs and arithmetic operations are in single precision and approximation code is accurate just enough to guarantee that casting of its double result to single precision is accurate enough. 
On average we observe that a slight performance improvement is still possible, sometimes even reaching 18\%. However, to achieve performance improvements comparable to those of double-precision code, we need a single-precision code generation from Metalibm.

{
	\renewcommand{\arraystretch}{1}
	\begin{table}[t]
		\centering
		\begin{tabular}{l|cc|cc|cc|cc|cc|cc}
			\toprule
			\multirow{3}{*}{benchmark} & \multicolumn{2}{c|}{\multirow{2}{*}{libm}} & \multicolumn{2}{c|}{small} & \multicolumn{2}{c|}{middle} & \multicolumn{6}{c}{large}                                                                 \\ \cline{4-13}
			& \multicolumn{2}{c|}{}                      & \multicolumn{2}{c|}{equal} & \multicolumn{2}{c|}{equal}  & \multicolumn{2}{c|}{equal} & \multicolumn{2}{c|}{derivative} & \multicolumn{2}{c}{no table} \\ \cline{2-13}
			& min                 & max                 & min         & max         & min          & max         & min         & max         & min            & max           & min           & max          \\ \hline
			axisRotationX     & 340 & 611 & 380 & 426 & 350 & 401 & 356 & 491 & 351 & 451 & 356 & 446 \\
      axisRotationY     & 337 & 541 & 378 & 480 & 352 & 473 & 352 & 419 & 294 & 509 & 352 & 428 \\
      ex2\_1            & 295 & 448 & 302 & 382 & 282 & 400 & 282 & 348 & 282 & 345 & 282 & 396 \\
      ex2\_10           & 339 & 510 & 296 & 338 & 324 & 354 & 324 & 339 & 326 & 448 & 324 & 342 \\
      ex2\_11           & 318 & 525 & 304 & 590 & 296 & 472 & 296 & 495 & 301 & 495 & 296 & 464 \\
      ex2\_2            & 299 & 490 & 373 & 415 & 340 & 368 & 294 & 363 & 294 & 348 & 294 & 368 \\
      ex2\_3            & 316 & 622 & 302 & 390 & 278 & 324 & 276 & 329 & 276 & 338 & 275 & 340 \\
      ex2\_4            & 296 & 493 & 310 & 332 & 298 & 349 & 288 & 356 & 291 & 371 & 288 & 345 \\
      ex2\_5            & 336 & 547 & 374 & 463 & 311 & 339 & 302 & 339 & 302 & 350 & 302 & 372 \\
      ex2\_9            & 320 & 561 & 358 & 382 & 346 & 390 & 346 & 445 & 346 & 401 & 346 & 402 \\
      ex3\_d            & 324 & 465 & 429 & 540 & 334 & 438 & 284 & 356 & 282 & 330 & 283 & 449 \\
      forwardk2jX       & 337 & 599 & 388 & 490 & 359 & 432 & 360 & 448 & 357 & 414 & 360 & 428 \\
      forwardk2jY       & 335 & 564 & 528 & 574 & 354 & 416 & 354 & 403 & 356 & 392 & 354 & 477 \\
      integrate18257    & 380 & 716 & 380 & 521 & 388 & 526 & 340 & 456 & 336 & 587 & 341 & 638 \\
      integStoutemyer   & 387 & 447 & 360 & 404 & 321 & 469 & 313 & 457 & 328 & 454 & 319 & 658 \\
      pendulum1         & 314 & 410 & 472 & 518 & 356 & 386 & 352 & 400 & 352 & 401 & 352 & 554 \\
      pendulum2         & 433 & 616 & 485 & 566 & 462 & 545 & 456 & 543 & 400 & 598 & 460 & 596 \\
      rodriguesRotation & 447 & 697 & 312 & 346 & 454 & 510 & 455 & 512 & 454 & 586 & 456 & 521 \\
      sinxx10           & 285 & 424 & 515 & 901 & 306 & 325 & 305 & 344 & 304 & 330 & 305 & 347 \\
      xu1               & 359 & 731 & 343 & 746 & 394 & 508 & 302 & 525 & 300 & 525 & 298 & 491 \\
      xu2               & 368 & 682 & 360 & 700 & 350 & 678 & 343 & 650 & 346 & 606 & 349 & 598 \\
    \bottomrule         
		\end{tabular}
		\caption{Minimum and maximum number of cycles.}\label{tab:minmax}
	\end{table}
	
}

\paragraph{Variance of Minimum and Maximum Cycles}
\autoref{tab:minmax} presents the average minimum and maximum number of cycles per program execution, in case of \libm-based implementation and of our new approximations. Our implementations in majority of cases provide significantly smaller variance between the minimum and maximum execution times. For large target errors the maximum execution time of our code is even often smaller than the best timings of \libm-based implementations. 
For instance, for the forwardk2jY benchmark, our the derivative-based is not only 17\% faster on average, but its maximum execution time is also considerably smaller than the one of the \libm-based code.

\begin{table}[t]
\scriptsize
\begin{tabular}{lrr|rrrrr}
\toprule
                        & small     & middle     & \multicolumn{5}{c}{large}                                           \\
benchmark               & equal     & equal      & equal      & deriv      & no table  & depth 1    & depth $\infty$  \\
\midrule
sinxx10                 & 30s    & 29s     & 32s     & 33s     &28s    & 7m 58s     & 6m 56s          \\
xu1                     & 6m 12s    & 4m 0s      & 10m 55s    & 8m 7s      & 10m 30s   & 10m 5s     & 8m 51s          \\
xu2                     & 15m 56s   & 7m 33s     & 3m 39s     & 3m 44s     & 3m 39s    & 3m 6s      & 4m 48s         \\
integrate18257          & 10m 1s    & 6m 49s     & 9m 1s      & 3m 24s     & 10m 19s   & 9m 1s      & 6m 52s         \\
integStoutemyer & 8m 23s    & 5m 25s     & 5m 19s     & 4m 18s     & 4m 56s    & 5m 39s     & 1m 10s          \\
axisRotationX           & 8m 58s    & 10m 42s    & 9m 13s     & 10m 45s    & 9m 11s    & -          & -               \\
axisRotationY           & 10m 46s   & 10m 44s    & 8m 22s     & 4m 2s      & 7m 43s    & -          & -               \\
rodriguesRotation       & 11m 5s    & 10m 43s    & 7m 28s     & 6m 14s     & 7m 30s    & 6m 57s     & 6m 59s           \\
pendulum1               &30s    &31s     &28s     &28s     &26s    & -          & -               \\
pendulum2               & 1m 22s    & 1m 24s     & 1m 18s     &54s     & 1m 57s    & 1m 36s     & 1m 30s          \\
forwardk2jX             & 5m 44s    & 3m 51s     & 7m 58s     & 4m 10s     & 6m 51s    & -          & -              \\
forwardk2jY             &56s    & 3m 22s     & 6m 4s      & 6m 18s     & 6m 2s     & -          & -               \\
ex2\_1                  & 5m 25s    & 12m 11s    & 5m 35s     & 5m 33s     & 5m 36s    & 5m 32s     & 5m 32s          \\
ex2\_2                  & 5m 32s    & 5m 31s     & 2m 52s     & 2m 35s     & 3m 28s    & 2m 25s     & 2m 25s           \\
ex2\_3                  & 10m 7s    & 7m 5s      & 8m 31s     & 8m 32s     & 8m 34s    & 8m 34s     & 8m 31s         \\
ex2\_4                  & 5m 36s    & 5m 54s     & 5m 33s     &50s     & 5m 33s    & 5m 33s     & 5m 34s           \\
ex2\_5                  &57s    & 1m 13s     &50s     &50s     &50s    &50s     &50s           \\
ex2\_9                  &55s    &50s     &49s     &48s     &48s    &49s     &49s           \\
ex2\_10                 &53s    &45s     &44s     &26s     &46s    &45s     &45s          \\
ex2\_11                 & 4m 51s    & 3m 7s      & 2m 5s      & 1m 39s     & 2m 4s     & 7m 1s      & 7m 4s           \\
ex3\_d                  & 5m 40s    & 5m 38s     & 4m 6s      & 4m 7s      & 4m 57s    & 4m 2s      & 4m 2s            \\
total                   & 2h19s & 1h 47m 47s & 1h 41m 22s & 1h 18m 17s & 1h 42m 6s & 1h 19m 54s & 1h 12m 38s  \\
\bottomrule   
\end{tabular}
\caption{Analysis time of our tool, measure using the bash time command.}
\label{tbl:analysis-time}
\end{table}

\paragraph{Analysis Time}
The time our tool takes for analysis and approximation is shown in~\autoref{tbl:analysis-time}.
Analysis time is highly dependent on the number of required approximations of elementary functions: each approximation requires a separate call to Metalibm whose running time, in its turn, depends on the problem definition. 
Daisy reduces the number of calls to Metalibm by common expression elimination which improves the analysis time.
Currently, we set the timeout for each Metalibm call to 3 minutes, which leads to an overall analysis time which is reasonable,
and at most 16min for our benchmarks.
We found it to be a reasonable bound for our relatively ``well-behaved'' functions. If Metalibm does not find an approximation within a few minutes, it usually means that it ``wobbles'' between domain splittings.

\section{Related Work}

\paragraph{Floating-point Analysis and Optimization}

  Several static analysis tools exist which bound roundoff errors of
  floating-point arithmetic computations, including elementary
  functions~\cite{real2Float,Precisa,FPTaylor}, but all assume libm library
  implementations and do not attempt to optimize the running time.
  Lee et. al.~\cite{Lee2017} verify the correctness of several existing \libm
  implementations in Intel's math library (exp, sin, tan, and log)
  automatically, but do not assume any approximations.
  
  Mixed-precision tuning~\cite{Precimonious,FPTuner,Lam2013b,DaisyTuning,Damouche2018}
  selects different floating-point precisions for different arithmetic
  operations such that an overall error bound is satisfied and performance is
  improved. This work has a similar premise as ours, however only considers the
  precision of arithmetic operations, and not of elementary functions. In
  practice, mixed-precision tuning is most helpful when an accuracy bound is
  desired which is close to what uniform precision can provide. In contrast, our
  approximation of elementary functions operates in a different part of the
  tradeoff space.

  The tool Herbie~\cite{Herbie} and Damouche et.al.~\cite{Damouche2017} perform a different kind of optimization.
  Instead of aiming to improve performance, it aims to increase the accuracy of
  a computation by leveraging the non-associativity of floating-point
  arithmetic. 

  Another way to improve the performance of (numerical) computations is
  autotuning, which performs low-level real-value semantics-preserving
  transformations of a program in order to find one which empirically executes
  most efficiently~\cite{Spiral,Vuduc2004}. The approaches optimize for
  different hardware platform and do not consider reducing accuracy in exchange
  for performance.

\paragraph{Elementary Function Approximation} 
The problem of evaluation of elementary functions is highly dependent on the technology. 
The theoretical support followed the needs of evolving hardware and software but was proceeding in a rather function-by-function manner. A good overview of the existing approaches for can be found in~\cite{Muller2016}.
However, to the best of our knowledge, Metalibm is the first and so far the only project aiming at automation of the \libm development and providing strong accuracy guarantees.

\paragraph{Approximate Computing}
  
  Approximate computing operates on the premise that many applications are
  tolerant to a certain amount of noise and errors and thus high accuracy is not
  always needed. Many techniques which trade accuracy for efficiency have been 
  developed in the recent past and are for instance surveyed in~\cite{ApproxSurvey}.

  Of particular interest to this project is work such as
  STOKE~\cite{StokeFloat}, which uses MCMC search to find reduced-precision
  implementations of short numerical kernels. Due to the stochastic search, the
  approach does not guarantee error bounds, however. The verifiable version of
  STOKE~\cite{STOKE}, which does not aim to reduce precision could be
  potentially used to improve the implementations generated by our tool further.

  The sketching synthesis technique has been used to, among others, synthesize
  polynomial approximations of programs with elementary
  functions~\cite{Metasketching}. Correctness, respectively accuracy, is only
  checked on a small set of sample inputs and thus cannot provide any
  guarantees.

  Similarly, neural networks have been used to learn approximations of numerical
  programs~\cite{NeuralApprox}, which can be efficiently executed on custom
  hardware, but again do not provide accuracy guarantees.

\section{Conclusion}

We presented a fully automated approach and tool which approximates elementary
functions inside small programs and provides rigorous whole-program error
guarantees which take into account both approximation as well as roundoff
errors. Our results show that it is possible to achieve significant performance
improvements in exchange for reduced, but guaranteed, accuracy.

Elementary function approximations are challenging and their efficient implementations
rely on a careful selection of many different parameters. Metalibm selects some
of these parameters, and our combination with Daisy's error analysis and
infrastructure allows our tool to select two more: correct error bounds of
individual elementary function calls, as well as the degrees of polynomial
approximations. This combination allows even non-experts to use rigorous
elementary function approximations in their programs.

While our tool already provides significant performance improvements, we note
that more work remains for the future: support for single-precision
implementations and multi-variate functions as well as more control over
Metalibm's heuristics, which will enable efficient implementations for even more complicated programs.

\section*{Acknowledgments}\label{sec:Acknowledgments}
The authors thank Christoph Lauter for useful discussions and Youcef Merah for the work on an early prototype.

\bibliographystyle{splncs04}
\bibliography{bibliography}
\newpage
\appendix
\section*{Appendix}

\paragraph{Benchmarks}
Following are the benchmarks that we used. The error bounds shown are for the
``small'' error case. The ``middle'' and ``large'' errors are obtained by increasing the bound by one and two orders of magnitude respectively.
\begin{lstlisting}
xu1: 
input: x1, x2 $\in$ [-3.14, 3.14]
  2 * sin(x1) + 0.8 * cos(2 * x1) + 7 * sin(x2) - x1
error bound:  1e-13 

xu2: 
input:  x1, x2 $\in$ [-3.14, 3.14]
  1.4 * sin(3 * x2) + 3.1 * cos(2 * x2) - x2 + 4 * sin(2 * x1)
error bound: 1e-13

integrate18257:
input x $\in$ [0, 3.14]
  exp(cos(x)) * cos(x - sin(x))
error bound: 1e-13

integrateStoutemyer2007:
input: x $\in$ [0.1, 1]
  log((exp(x) + 2 * sqrt(x) + 1) / 2.0)
error bound: 1e-13

sinxx10:
input: x $\in$ [-3, 3]
  (3 * x * x * x - 5 * x + 2) * sin(x) * sin(x) + (x * x * x + 5 * x) * sin(x) - 2*x*x - x - 2
error bound: 1e-11

axisRotationX:
input: x $\in$ [-2, 2], y $\in$ [-4, 4], theta $\in$ [-5, 5]
  x * cos(theta) + y * sin(theta)
error bound: 1e-13

axisRotationY:
input: x $\in$ [-2, 2], y $\in$ [-10, 10], theta $\in$ [-5, 5]
  -x * sin(theta) + y * cos(theta)
error bound: 1e-12

rodriguesRotation: 
input: (v1 $\in$ [-2, 2], v2 $\in$ [-2, 2], v3 $\in$ [-2, 2], k1 $\in$ [-5, 5], k2 $\in$ [-5, 5], 
k3 $\in$ [-5,5], theta $\in$ [-5, 5]
  v1 * cos(theta) + (k2 * v3 - k3 * v2) * sin(theta) + 
    k1 * (k1 * v1 + k2 * v2 + k3 * v3) * (1 - cos(theta))
error bound: 1e-11

pendulum1:
input: t $\in$ [1, 3], w $\in$ [-5, 5]
  t + 0.01 * (w + 0.01/2*(-9.80665/2.0 * sin(t)))
error bound: 1e-13

pendulum1:
input: t $\in$ [-2, 2], w $\in$ [1, 5]
  w + 0.01 * exp(-9.80665/2.0 * sin(t + 0.01/2*cos(w)))
error bound: 1e-13

forwardk2jY:
input: theta1 $\in$ [-3.14, 3.14], theta2 $\in$ [-3.14, 3.14]
  0.5 * sin(theta1) + 2.5 * sin(theta1 + theta2)
error bound: 1e-13


forwardk2jX:
input: theta1 $\in$ [-3.14, 3.14], theta2 $\in$ [-3.14, 3.14]
  0.5 * cos(theta1) + 5.5 * cos(theta1 + theta2)
error bound: 1e-13

ex2_1: 
input: x $\in$ [-1.57, 1.57]
  val x1 = sin(x); val x2 = cos(x); x1 * x1 * x2 * x2
error bound: 1e-13

ex2_2:
input: x $\in$ [-1, 1]
  val x1 = sin(2 * x);  val x2 = cos(2 * x);  x1 * x1 * x2 * x2 * x2
error bounnd: 1e-13

ex2_3:
input: x $\in$ [0, 1]
  val x1 = cos(2 * x); val x2 = cos(3 * x); x1 * x2
error bound: 1e-13

ex2_4:
input: x $\in$ [-1.57, 1.57]
  val x1 = sin(x); val x2 = cos(x); x1 * x1 * x1 * x1 * x1 * x2 * x2
error bound: 1e-13

ex2_5: 
input: x $\in$ [17, 18] 
  val x1 = sin(x); val x2 = cos(x); (x1 + 2 * x2) / (x2 + 2 * x1)
error bound: 1e-12

ex2_9:
input: x $\in$ [1, 3.1415]
  val x1 = sin(x); val x2 = cos(x); 1 / (1 - x2 + x1)
error bound: 1e-13

ex2_10:
input: x $\in$ [-20, -18]
  val x1 = sin(x); val x2 = 1 + cos(x); x1 / (x2 * x2)
error bound: 1e-13

ex2_11:
input: x $\in$ [-1.1, 0.9]
  val x1 = 1 / cos(x); val x2 = tan(x); (x1 * x1) / (4 + x2 * x2)
error bound: 1e-12

ex3_d
input: x $\in$ [0, 7]
  val x1 = exp(-2 * x); val x2 = sin(x); x1 * x2
error bound: 1e-13
\end{lstlisting}

\end{document}